\documentclass[12pt]{amsart}
\usepackage{amssymb}
\usepackage{amsthm}
\usepackage{amsmath}
\usepackage{enumerate}
\usepackage{cancel}
\usepackage[normalem]{ulem}
\usepackage[utf8]{inputenc}

\newtheorem{theorem}{Theorem}
\newtheorem{remark}{Remark}

\newtheorem{corollary}{Corollary}
\newtheorem{example}{Example}

\begin{document}
	\title[Bases of the Riemann-Roch spaces on divisors
	on hyperelliptic curves]
	{Explicit bases of the Riemann-Roch spaces on divisors
	on hyperelliptic curves}
	
	\thanks{The first author was supported by University of Palermo FFR. The second author was supported by the National Research, Development and Innovation Office (NKFIH) Grant No. K132951 and by the contruction EFOP-3.6.1-16-2016-00022, which is a project supported by the European Union, co-financed by the European Social Fund. 
	The third author was supported by the construction EFOP-3.6.3-VEKOP-16-2017-00002. These projects were supported by the European Union, co-financed by the European Social Fund. }
	
	\author[G. Falcone]{Giovanni Falcone} \address{Department of Mathematics and Computer science, University of Palermo, Italy} \email{giovanni.falcone@unipa.it}
	\email{ORCID: 0000-0002-5210-5416}
	
	\author[\'A. Figula]{\'Agota Figula}
	\address{Institute of Mathematics, University of Debrecen, Hungary} \email{figula@science.unideb.hu}
	\email{ORCID: 0000-0002-8095-6074}
	
	\author[C. Hannusch]{Carolin Hannusch}
	\address{Faculty of Informatics, University of Debrecen, Hungary}
	\email{hannusch.carolin@inf.unideb.hu}
	\email{ORCID: 0000-0002-0098-7293}

	
	
	\date{}

\begin{abstract} 
For an (imaginary) hyperelliptic curve $\mathcal{H}$ of genus $g$, we determine a basis of the Riemann-Roch space
$\mathcal{L}(D)$, where
$D$ is a divisor with positive degree $n$, linearly equivalent to $P_1+\cdots+ P_j+(n-j)\Omega$, with
$0 \le j \le g$,
where $\Omega$ is a Weierstrass point, taken as the point at infinity. As an application, we determine a generator matrix of a Goppa code for $j=g=3$ and $n=4.$	
\end{abstract}



\maketitle

{\small\centerline{In memoriam Prof. Dr. Heinrich Wefelscheid (16.4.1941-18.4.2020)}}

\bigskip

\section{Introduction}
\label{introduction}
Let $\chi $ be an algebraic curve defined over a finite field $\mathrm{GF}(q)$,  and let $D$ be an
$\mathrm{GF}(q)$-rational divisor on $\chi $.
The computation of a basis of the Riemann-Roch space $\mathcal{L}(D)$ associated to $D$ is an essential tool in Coding Theory and Cryptography, since it allows both to explicitly construct Goppa codes and to give addition formulas in the divisor class group of $\chi $. The general problem has been attacked by several researchers. The first algorithm is due to von Brill and Noether \cite{BrillNoether}. Since then, many researchers have worked on the problem to make the computation of such a basis more effective (\cite{HuangIerardi}), in the equivalent scenario of function fields 
(cf. \cite{Stichtenoth}, Remark 2.3.15). In particular, in \cite{Hess} an arithmetic
approach to the Riemann-Roch problem is taken, which provides an algorithm, polynomial in the input size. 
Nonetheless, further algorithms were developed in order to simplify the computation, each under particular assumptions.

In this paper the class of hyperelliptic curves is considered. Many papers have been devoted to the study of arithmetic in these curves, among others we mention in particular \cite{Cantor}, \cite{kuroki} and \cite{Lange2}. The interest on the subject does not seem to decline, as witnessed by more recent publications (cf. \cite{sutherland}, \cite{GluherSpaenlehauer}). A significant literature has also been produced in order to consider codes over hyperelliptic curves (\cite{boer0}, 
\cite{boer}, \cite{brigand},  \cite{niehage2}, 

Goppa codes were introduced in \cite{G2} several decades ago. These codes turned out not only to be interesting in Coding Theory, but also to be applicable in Cryptography, e.g. in public-key cryptographic systems \cite{Mc}, \cite{JanwaMoreno}, \cite{MMPR}. Cryptographic systems using Goppa codes with suitable parameters are considered as secure (see in particular \cite{BLP}, \cite{DinhMooreRussell}). Hyperelliptic curves in Cryptography have been  investigated in \cite{koblitz}, \cite{kuroki}, 
\cite{pelzl}. Goppa codes over the Hermitian curve, as well as over maximal curves, have been extendedly studied in \cite{KorchmarosNagyTimpanella}, \cite{KorchmarosSpeziali}, \cite{CastellanosFanali}, \cite{FanaliGiulietti}.

The aim of the current paper is to give an explicit way for determining a basis of the Riemann-Roch space over an imaginary hyperelliptic curve $\mathcal{H}.$ Using this basis, we can construct a generator matrix of a Goppa code over a hyperelliptic curve defined over a Galois field of characteristic $p\geq 2$. We make this for some MDS codes in Section \ref{ExampleGoppa}. In particular we consider an imaginary hyperelliptic curve
$\mathcal{H}$ of genus $g$, described as the set of points satisfying the equation
$$Y^2T^{d-2}+Yh(X,T)=f(X,T)$$
where $f$ is a homogeneous polynomial of degree $d=2g+1$ (and $h=0$, if $p \neq 2$). Using standard methods we construct an explicit  basis of the Riemann-Roch space $\mathcal{L}(D)$, where $D$ is a divisor of positive degree 
$n$ in (its unique) reduced form $P_1+\dots +P_j+(n-j)\Omega$. Here $P_1, \dots , P_j$ are $j$ points in 
$\mathcal{H}$ distinct from the point $\Omega$ of infinity and $j \le g$.  We remark that the reduction of $D$ to its reduced form might be a difficult task, because one has to solve algebraic equations of degree greater than 
$g+1$, possibly by applying the Cantor algorithm. This difficulty does not occur in the construction of Goppa codes, because in that case one can directly take $D=P_1+P_2+\dots +P_j+(n-j)\Omega$ (cf. \cite{Cantor}).

It turns out that our computation for a basis of the Riemann-Roch space $\mathcal{L}(D)$ provides another proof of the results of Lemma 2.1 in \cite{boer}, which deal with the dimension of the space $\mathcal{L}(D)$. We give the sequence of 
$\mathrm{dim} \mathcal{L}(D)$ in Example 1 for the case that $\mathcal{H}$ has genus $g=5$.

\section{Notations and definitions}
\label{sec:1}
Let $p$ be a prime number and $t\in \mathbb{N}.$
Let $\mathcal{H}$ be a hyperelliptic curve over $\mathrm{GF}(p^t)$ with a rational Weierstrass point $\Omega$, so that there exists a coordinate system of the projective plane such that the non-singular curve $\mathcal{H}$ is described as the set of points $P=[X:Y:T]$ such that
$$Y^2T^{d-2}+Yh(X,T)=f(X,T)$$
where $f$ is a homogeneous polynomial of degree $d=2g+1$, $h$ is a homogeneous polynomial of degree at most $g$, and
$\Omega=[0:1:0]$ is the point at infinity of $\mathcal{H}$ (\cite{Lockhart}, Prop. 1.2). If $p$ is odd, the transformation
$Y\mapsto Y-h(X,T)/2$ changes the above equation into
$$Y^2T^{d-2}=f(X,T),$$
whereas, if $p=2$, then in general it is not possible to reduce $h$ to zero.

Let $\mathfrak{K}$ be the algebraic closure of $\mathrm{GF}(p^t)$, and let $\mathcal{L}(D)$ be the Riemann-Roch space associated to any divisor $D$, that is, the vector space of rational functions
$$\mathcal{L}(D)=\{F\in \mathfrak{K}(\mathcal{H}):\mathrm{div}(F)+D\geq 0\}\cup\{0\},$$
thus $\mathcal{L}(D)$ is trivial both in the cases where $D$ has negative degree, and where $D$ has degree zero and
$D\not\in\mathrm{Princ}(\mathcal{H})$, whereas $\mathcal{L}(D)=\Big\langle  F_0^{-1} \Big\rangle  $ in the case where
$D=\mathrm{div}(F_0)$.
For this reason, we may restrict ourselves to the case where $D$ has positive degree.

If $D$ is a divisor of positive degree $n$, then
$$D=P_1+P_2+\dots +P_j+(n-j)\Omega+\mathrm{div}(\psi)$$
for $j$ points $P_1,\dots , P_j$ in $\mathcal{H}$ distinct from $\Omega$, with $j\leq g$, and a suitable 
$\psi\in{\mathfrak{K}}(\mathcal{H})$, that is,
any divisor class $D+\mathrm{Princ}(\mathcal{H})\in \mathrm{Div}(\mathcal{H})/\mathrm{Princ}(\mathcal{H})$ can be reduced to the form $P_1+\dots +P_j+(n-j)\Omega$.

Up to the isomorphism $$\Phi: \mathcal{L}(D)\mapsto \mathcal{L}(P_1+\dots +P_j+(n-j)\Omega),$$ mapping $F$ onto the product
$\psi F$, we will directly assume that $D$ is
reduced to $P_1+\dots +P_j+(n-j)\Omega$, $n \ge 0$.

\section{Main theorem}\label{Main} 
Let $\mathcal{H}$ be the hyperelliptic curve introduced in Section \ref{sec:1}. Let $D=P_1+\dots +P_j+(n-j)\Omega$ be a divisor of degree $n$ of $\mathcal{H}$.

If $P_i=[a_i:b_i:1]$, then let $Q_i=[a_i:-b_i-h(a_i,1):1]\in\mathcal{H}$, and let
\begin{equation}\label{Psi}
\Psi=\frac{T^{j-\delta}\kappa}{(X-a_1T)\cdots(X-a_jT)},
\end{equation}
where $\kappa$ is the curve $YT^{\delta-1}-k(X,T)$ of smallest degree $\delta$ in $X$ passing through the points
$Q_1,\dots,Q_j$ with their possible multiplicity (note, in particular, that for $j=1$ the curve $\kappa$ is the line $Y-(b_1+h(a_1,1))T$, and recall that $h(X,T)=0$, if $p>2$). Furthermore, we define $\Psi=\frac{Y}{T}$ for 
$j=0.$

Since $\delta<j\leq g=\frac{d-1}{2}$, there are $d$ intersection points of $\kappa$ and $\mathcal{H}$ in the affine plane, say $Q_1,\dots, Q_j$ and $W_1,\dots,W_{2g-j+1}$ (in the case where $j=0$, these being the $d$ intersections of $\mathcal{H}$ with the $x$-axis), and $d\cdot (\delta -1)$ further intersection points coinciding with $\Omega$, hence
$$\mathrm{div}\Psi=(W_1+\dots +W_{2g-j+1})-(P_1+\dots+P_j)-(2g-2j+1)\Omega.$$

\begin{remark}\label{Rmkpsi}
\emph{	Note that $\Psi\in\mathcal{L}(D)$ if and only if $n-j\geq 2(g-j)+1$.}
\end{remark}

\begin{theorem}\label{MainTheorem} Let $D=P_1+\dots+P_j+(n-j)\Omega$ be a divisor of degree $n$ on the hyperelliptic curve
	$\mathcal{H}$ defined in Section \ref{sec:1}, and let $\Psi$ be as in (\ref{Psi}).
	If $n-j\geq 2(g-j)+1$, then a basis of $\mathcal{L}(D)$ is provided by the set
	$$\left\{ \left(\frac{X}{T}\right)^h,\Psi\left(\frac{X}{T}\right)^k:\;
	0\leq h\leq\frac{n-j}{2}\mbox{ and }0\leq k\leq\frac{(n-j)-2(g-j)-1}{2}\right\}.$$
	If $n-j<2(g-j)+1$, then a basis of $\mathcal{L}(D)$ is provided by the set
	$$\left\{ \left(\frac{X}{T}\right)^h:\;
	0\leq h\leq\frac{n-j}{2}\right\}.$$
\end{theorem}

\begin{proof}
	Let $B_1$ and $B_2$ be the intersection points of $\mathcal{H}$ and the $y$-axis, so
	$$\mathrm{div}\left(\frac{X}{T}\right)=(B_1+B_2)-2\Omega.$$
	
	1) Let $n-j\geq 2(g-j)+1$, thus $\Psi\in\mathcal{L}(D)$.
	First we consider the cases where either $j=0$ (hence $n \ge 2g+1$), or  $j=1$ (hence $n \ge 2g$), or $j \ge 2$ and $n\geq 2g-1$, as in these cases we know that, by the theorem of Riemann-Roch, the dimension of
	$\mathcal{L}(D)$ is $n-g+1$. We claim that
	$$\mathcal{L}(D)=\left\langle \left(\frac{X}{T}\right)^h,\Psi\left(\frac{X}{T}\right)^k\right\rangle,\mbox{ where}$$
	$$0\leq h\leq\frac{n-j}{2}\mbox{ and }0\leq k\leq\frac{(n-j)-2(g-j)-1}{2}.$$
	In fact, for each of those values of the parameters $h$ and $k$, the functions belong to $\mathcal{L}(D)$, because
	$$D+\mathrm{div}\left(\frac{X}{T}\right)^h=(P_1+\dots+P_j)+(n-j)\Omega+h(B_1+B_2)-2h\Omega,$$
	as well as
	$$D+\mathrm{div}\;\Psi\left(\frac{X}{T}\right)^k=(P_1+\dots+P_j)+(n-j)\Omega+k(B_1+B_2)-2k\Omega +$$
	$$+(W_1+\dots +W_{2g-j+1})-(P_1+\dots+P_j)-(2g-2j+1)\Omega$$
	$$=k(B_1+B_2)+(W_1+\dots +W_{2g-j+1})-(2k-(n-j)+2(g-j)+1)\Omega,$$
	are effective divisors.
	Since $0\leq h\leq\frac{n-j}{2}$ and $0\leq k\leq\frac{(n-j)-2(g-j)-1}{2}$, if $(n-j)$ is even, then the number of those functions is
	$$1+\frac{n-j}{2}+1+\frac{(n-j)-2(g-j)-2}{2}=n-g+1,$$
	and, if $(n-j)$ is odd, then their number is
	$$1+\frac{n-j-1}{2}+1+\frac{(n-j)-2(g-j)-1}{2}=n-g+1,$$
	as well, and the claim follows from dimensional reasons.
	
	\medskip
	Secondly, we consider the case where $2g-j+1\leq n< 2g-1$ (note that this case can occur only if $j\geq 3$). In this case, the dimension of $\mathcal{L}(D)$ is not necessarily $n-g+1$, but still $\Psi\in\mathcal{L}(D)$, and again we claim that
	$$\mathcal{L}(D)=\left\langle \left(\frac{X}{T}\right)^h,\Psi\left(\frac{X}{T}\right)^k\right\rangle,$$
	where $0\leq h\leq\frac{n-j}{2}$ and $0\leq k\leq\frac{(n-j)-2(g-j)-1}{2}$.
	
	In fact, let $n=2g-1-\epsilon$ with $0\leq\epsilon\leq j-2$, and put, for short,
	$$\mathcal{L}_\epsilon:=\mathcal{L}(P_1+\dots+P_j+(n-j)\Omega),$$
	hence we have from the first case above that
	$$\mathcal{L}_0=\left\langle \left(\frac{X}{T}\right)^h,\Psi\left(\frac{X}{T}\right)^k\right\rangle$$
	where $0\leq h\leq\frac{n-j}{2}$ and $0\leq k\leq\frac{(n-j)-2(g-j)-1}{2}$. Since
	$$\mathcal{L}_{\epsilon+1}\leq \mathcal{L}_\epsilon,$$
	it follows recursively, by step by step inspection of $\mathcal{L}_\epsilon$, that the claim holds as in the first case.
	
	\bigskip\noindent
	2) Let $j=0$ and $n=2g-1, 2g$, or $j=1$ and $n=2g-1$ (hence, by Remark \ref{Rmkpsi}, $\Psi\not\in\mathcal{L}(D)$).
	Again by the theorem of Riemann-Roch, the dimension of $\mathcal{L}(D)$ is $n-g+1$ and we have, by dimensional reason, that
	$\mathcal{L}(D)=\left\langle \left(\frac{X}{T}\right)^h\right\rangle$, where $0\leq h\leq\frac{n}{2}$, respectively $0\leq h\leq\frac{n-1}{2}$.

	\bigskip\noindent
	3) Finally, let either $j=0, 1$ and $n <2g-1$, or $2\leq j\leq n<2g-j+1$. In either of these cases, we claim that
	$\mathcal{L}(D)=\left\langle \left(\frac{X}{T}\right)^h\right\rangle$, where $0\leq h\leq\frac{(n-j)}{2}$. Let
	$n=2g-j+1-\epsilon$ with $0\leq\epsilon\leq 2g-2j+1$, and again put, for short,
	$$\mathcal{L}_\epsilon:=\mathcal{L}(P_1+\dots+P_j+(n-j)\Omega),$$
	hence we have from the first two cases above that
	$$\mathcal{L}_0=\left\langle \left(\frac{X}{T}\right)^h,\Psi\right\rangle$$
	where $0\leq h\leq\frac{n-j}{2}$, because for $n=2g-j+1$ we get $k=0$. Since, again,
	$$\mathcal{L}_{\epsilon+1}\leq \mathcal{L}_\epsilon,$$
	and since, by Remark \ref{Rmkpsi}, $\Psi\not\in\mathcal{L}_\epsilon$ as soon as $\epsilon>0$,
	it follows recursively, by step by step inspection of $\mathcal{L}_\epsilon$, that the claim holds as in the first case.\end{proof}

\begin{corollary} \label{dim}
	Let $\mathcal{H}$ be the hyperelliptic curve of genus $g$ defined in Section \ref{sec:1}. If the divisor $D$ of degree $n$ is linearly equivalent to $P_1+\dots+P_j+(n-j)\Omega$, then
	
	$\mathrm{dim}\, \mathcal{L}(D)=n-g+1$, for $n \geq 2g-j$,
	
	$\mathrm{dim}\, \mathcal{L}(D)=\lfloor\frac{n-j}{2}\rfloor +1$, for $j\leq n < 2g-j$.
	
	\noindent
	In particular, if $j=g$, then $\mathrm{dim}\, \mathcal{L}(D)=n-g+1$.
\end{corollary}

The above results on the dimension of the Riemann-Roch space $\mathcal{L}(D)$ appeared first in \cite{boer}, Lemma 2.1.

\begin{remark}
\emph{	A point $P$ of the curve $\mathcal{H}$ is a non-Weierstrass point if the sequence $\mathrm{dim}\mathcal{L}(nP)$ for $n \ge 1$ is
		$$\underbrace{1,1,\dots ,1}_g,2,3,4,\dots ,g-1,g,g+1,\dots$$
		In any other case, $P$ is a Weierstrass point. 
		Recall that $\Omega$ was a given Weierstrass point of the curve $\mathcal{H}$, and in fact
		$\mathcal{L}\big((2g-2\epsilon)\Omega\big)=\mathcal{L}\big((2(g-\epsilon)+1)\Omega\big)$ have both dimension
		$g-\epsilon+1$. The sequence $\mathrm{dim}\,\mathcal{L}(n \Omega)$, where $n \ge 0$, is therefore:
		$$1,1,2,2,3,3,\dots, g-1,g-1, g, g, g+1,g+2,g+3\dots$$
		Thus the sequence of gaps one has to fill from the first entry to any increasing entry is
		$$1,3,5,\dots, 2g-1,$$
		and the numerical semigroup of non-gaps is therefore that of the natural numbers without the odd numbers smaller than $2g$. }
\end{remark}

\begin{example} \emph{ Assume that $\mathcal{H}$ has genus $g=5$. }

\emph{ If $j=0$, then $n-j \ge 2(g-j)+1$ if and only if $n \ge 2g+1$. Hence the sequence of $\mathrm{dim}\mathcal{L}(D)$ is
$$\underbrace{1,1,2,2,3,3,4,4,5,}_{0\leq n<2g-1 } \underbrace{5,6,7, \dots}_{n \ge 2g-1}$$ }

\emph{ If $j=1$, then $n-j \ge 2(g-j)+1$ if and only if $n \ge 2g$. Hence the sequence of $\mathrm{dim}\mathcal{L}(D)$ is
$$\underbrace{1,1,2,2,3,3,4,4,}_{1\leq n<2g-1 } \underbrace{5,6,7, \dots}_{n \ge 2g-1}$$ }

\emph{ If $j=2$, then $n-j \ge 2(g-j)+1$ if and only if $n \ge 2g-1$. Hence the sequence of $\mathrm{dim}\mathcal{L}(D)$ is
$$\underbrace{1,1,2,2,3,3,4,}_{2\leq n<2g-1 } \underbrace{5,6,7, \dots}_{n \ge 2g-1}$$ }

\emph{ If $j=3$, then $n-j \ge 2(g-j)+1$ if and only if $n \ge 2g-2$. Hence the sequence of $\mathrm{dim}\mathcal{L}(D)$ is
$$\underbrace{1,1,2,2,3,}_{3\leq n<2g-2 }{4,} \underbrace{5,6,7, \dots}_{n\geq 2g-1}$$ }

\emph{ If $j=4$, then $n-j \ge 2(g-j)+1$ if and only if $n \ge 2g-3$. Hence the sequence of $\mathrm{dim}\mathcal{L}(D)$ is
$$\underbrace{1,1,2}_{4\leq n<2g-j+1 }\; \underbrace{3,4}_{2g-j+1\leq n<2g-1 }\; \underbrace{5,6,7, \dots}_{n\geq 2g-1}$$ }

\emph{If $j=5$, then $n-j \ge 2(g-j)+1$ if and only if $n \ge 2g-4$. Hence the sequence of $\mathrm{dim}\mathcal{L}(D)$ is
$$\underbrace{1,}_{5\leq n<2g-j+1 } \;\underbrace{2,3,4}_{2g-j+1\leq n<2g-1 }\; \underbrace{5,6,7, \dots}_{n\geq 2g-1}$$ } \end{example}

\begin{remark}
\emph{	Note that we do not need to know the equation of $\mathcal{H}$ of genus $g$ in order to get the basis of $\mathcal{L}(D).$}
\end{remark}

\section{Examples of MDS Goppa codes}\label{ExampleGoppa}
We take $j=3$, therefore $g=3$ and $n= 4$ for an example where $2g-j+1\leq n<2g-1$. Choose, for instance, $D=[0:0:1]+2[1:0:1]+\Omega$ over the Galois field $\mathrm{GF}(31)$, where we intentionally took the point $[1:0:1]$ twice. With the notation in the proof of Theorem \ref{MainTheorem}, put
$$\begin{array}{lll}
Q_1\equiv[0:0:1]; & Q_2\equiv[1:0:1]; & Q_3=Q_2. 
\end{array}$$
As $Q_3=Q_2$, we give the parabola $\kappa$ in the form
$$(y-0)=a_1(x-1)+a_2(x-1)^2,$$
and because it passes through $Q_1$, we obtain $a_2=a_1$, so take $y=(x-1)+(x-1)^2$. Putting $\Psi=\frac{YT^2-T\big(T(X-T)+(X-T)^2\big)}{X(X-T)^2}$, from Theorem \ref{MainTheorem} we obtain
$$\mathcal{L}(D)=\left\langle 1,\Psi\right\rangle.$$
We want to construct the $(m,2,\mathfrak{d})$-Goppa code (where $ m-4\leq\mathfrak{d}\leq m-1$)\footnote{The minimal distance of a Goppa code is $\mathfrak{d}\geq m-\mathrm{deg}(D)$.}. Therefore we need the equation of $\mathcal{H},$ so that we can take $m$ points on $\mathcal{H}.$ Thus we choose four further points $W_1,\dots, W_4$ on the parabola $\kappa$, and a fifth point $C$ not belonging to $\kappa$.
If we take, for instance, the abscissa of $W_i$ from $3$ to $6$, then we obtain the following points on $\kappa$:
$$\begin{array}{ll}
W_1\equiv[3:6:1] & W_2\equiv[4:12:1]  \\
W_3\equiv[5:20:1] & W_4\equiv[6:30:1]
\end{array}.$$
Since the point on the parabola with abscissa $x=7$ is $[7:11:1]$, we take $C=[7:12:1]$.

Thereafter, in order to have the equation of $\mathcal{H}$ in the form $$y^2=a_2(x-1)^2+\dots +a_7(x-1)^7,$$
we construct the $6\times 6$ Vandermonde matrix $V=\big((x_i-1)^{1+j}\big)$ of the absciss{\ae} $x_i$ of the points $Q_1,W_1,\dots,W_4,C$, and its inverse, that is,
{\small{$$V=\left(
		\begin{array}{cccccc}
		1&30&1&30&1&30\\
		4&8&16&1&2&4\\
		9&27&19&26&16&17\\
		16&2&8&1&4&16\\
		25&1&5&25&1&5\\
		5&30&25&26&1&6\end{array}
		\right),\;V^{-1}=\left(
		\begin{array}{cccccc}
		18&9&23&18&11&9\\  8&2&3&15&5&30\\  30&20&22&3&7&24\\  0&
		25&25&23&28&19\\  16&5&15&14&14&6\\  24&7&1&28&30&21
		\end{array}
		\right).$$}} 
Since $V^{-1}[0^2,6^2,12^2,20^2,30^2,12^2]'= [22, 10, 26, 3, 14, 18]'$, the hyperelliptic curve $\mathcal{H}$ is defined by the equation
$$y^2=22(x-1)^2+10(x-1)^3+26(x-1)^4+3(x-1)^5+14(x-1)^6+ 18(x-1)^7.$$

\noindent
With respect to the divisor $G=[3:25:1]+[4:19:1]+[5:11:1]+[6:1:1] \in \mathcal{H}$, not in the support of $D$, the generator matrix of the
$(4,2)_{31}$-Goppa code is
$\mathcal{C}_{\mathcal{L}}(D,G)=\left( \begin{array}{cccc}
1&1&1&1\\
30&20&15&12 \end{array} \right)$. Hence, the parity-check matrix is $H=\left( \begin{array}{cccc}
16&14&1&0\\ 7&23&0&1 \end{array} \right)$, and one sees that two columns in $H$ are always independent, thus the minimum distance $\mathfrak{d}$ is $3,$ i.e. the code is MDS.

\begin{remark}
\emph{	Note that, for $p\leq\frac{n-j}{2}$, the polynomials $X/T$ and $(X/T)^p$ take the same values in the field $\mathrm{GF}(p)$, thus one does not obtain linearly independent row vectors in the generator matrix. }
\end{remark}

Now we compute the generator matrices of some MDS Goppa codes of dimension $3$ arising from hyperelliptic curves of genus $g=2$ and constructing by de Boer in \cite{boer1}, Section 2. 

\begin{example} \emph{According to Example $2.3.21$ in \cite{boer1} we take the hyperelliptic curve 
$\mathcal{H}: Y^2T^3=X^5+4 X T^4+T^5$ over the field $GF(5)$ and we choose the divisor $D=[0:1:1]+[1:4:1]+2\Omega$. With the notation in Section \ref{Main} we put
$$\begin{array}{lll}
Q_1\equiv[0:4:1]; & Q_2\equiv[1:1:1].  
\end{array}$$
The line $\kappa$ passing through $Q_1$ and $Q_2$ has the form
$y=2x+4$. We have $j=2$ and $n=4$. Putting $\Psi=\frac{T\big(Y+3X+T\big)}{X(X-T)}$ from Theorem \ref{MainTheorem} we obtain
$$\mathcal{L}(D)=\left\langle 1, \frac{X}{T}, \Psi\right\rangle.$$
With respect to the divisor $G=[2:1:1]+[2:4:1]+[3:1:1]+[3:4:1]+[4:1:1]+[4:4:1] \in \mathcal{H}$, not in the support of $D$, the generator matrix of the
$(6,3,4)_{5}$-$MDS$ code is
$\mathcal{C}_{\mathcal{L}}(D,G)=\left( \begin{array}{cccccc}
1&1&1&1&1&1\\
2&2&3&3&4&4\\
4&3&1&4&2&1 \end{array} \right)$.} \end{example} 

\begin{example} \emph{According to Example $2.3.23$ in \cite{boer1} we take the hyperelliptic curve 
$\mathcal{H}: Y^2T^3=X^5+4X^3T^2+9XT^4$ over the field $GF(13)$ and we choose the divisior $D=[0:0:1]+3\Omega$. With the notation in Section \ref{Main} we have $Q\equiv[0:0:1]$  
and $\kappa$ has the form $Y$ because of $j=1$.  Putting $\Psi=\frac{Y}{X}$ and taking into account that $n=4$    Theorem \ref{MainTheorem} yields that 
$\mathcal{L}(D)=\left\langle 1, \frac{X}{T}, \Psi\right\rangle.$
With respect to the divisor $G=[1:1:1]+[1:12:1]+[3:1:1]+[3:12:1]+[6:6:1]+[6:7:1]+[7:4:1]+[7:9:1]+[9:6:1]+[9:7:1] \in \mathcal{H}$, not in the support of $D$, the generator matrix of the
$(10,3,8)_{13}$-$MDS$ code is
$\mathcal{C}_{\mathcal{L}}(D,G)=\left( \begin{array}{cccccccccc}
1&1&1&1&1&1&1&1&1&1\\
1&1&3&3&6&6&7&7&9&9\\  
1&12&9&4&1&12&8&5&5&8 \end{array} \right)$.} \end{example}

\begin{example} \emph{According to Example $2.3.24$ in \cite{boer1} we take the hyperelliptic curve 
$\mathcal{H}: Y^2T^3=X^5+13X^4T+5X^3T^2+11X^2T^3+5XT^4+15T^5$ over the field $GF(17)$ and we choose the divisior $D=[8:0:1]+3\Omega$. With the notation in Section \ref{Main} we have $Q\equiv[8:0:1]$  
and as $j=1$ the form of $\kappa$ is again $Y$.  Putting $\Psi=\frac{Y}{X-8}$ from Theorem \ref{MainTheorem} we obtain
$\mathcal{L}(D)=\left\langle 1, \frac{X}{T}, \Psi\right\rangle $ since $n=4$. 
With respect to the divisor $G=[0:7:1]+[0:10:1]+[1:4:1]+[1:13:1]+[3:8:1]+[3:9:1]+[5:1:1]+[5:16:1]+[9:1:1]+[9:16:1]+ [15:7:1]+[15:10:1]\in \mathcal{H}$, not in the support of $D$, the generator matrix of the
$(12,3,10)_{17}$-$MDS$ code is
$\mathcal{C}_{\mathcal{L}}(D,G)=\left( \begin{array}{cccccccccccc}
1&1&1&1&1&1&1&1&1&1&1&1\\
0&0&1&1&3&3&5&5&9&9&15&15\\
14&3&14&3&12&5&11&6&1&16&1&16 \end{array} \right)$.} \end{example} 

\begin{example} \emph{According to Examples $2.2.3$, $2.3.26$ in \cite{boer1} we take the hyperelliptic curve 
$\mathcal{H}: Y^2T^3+YT^4=X^5+X^3T^2+XT^4$ over the field $GF(4)$ and we choose the divisior $D=[\alpha^2:0:1]+3\Omega$, where 
$\alpha $ is a primitive element satisfying the equation $\alpha^2+\alpha+1=0$. With the notation in Section \ref{Main} we have 
$Q\equiv[\alpha^2:1:1]$  
and as $j=1$ the form of $\kappa$ is $Y-T$.  Putting $\Psi=\frac{Y-T}{X-\alpha^2T}$ from Theorem \ref{MainTheorem} we obtain
$\mathcal{L}(D)=\left\langle 1, \frac{X}{T}, \Psi\right\rangle $ since $n=4$. 
With respect to the divisor $G=[0:0:1]+[0:1:1]+[1:\alpha:1]+[1:\alpha^2:1]+[\alpha:0:1]+[\alpha:1:1] \in \mathcal{H}$, not in the support of $D$, the generator matrix of the
$(6,3,4)_{4}$ hexacode is
$\mathcal{C}_{\mathcal{L}}(D,G)=\left( \begin{array}{cccccc}
1&1&1&1&1&1\\
0&0&1&1&\alpha&\alpha\\
\alpha&0&\alpha&1&1&0 \end{array} \right)$.} \end{example}


\begin{thebibliography}{ABC}
		
	\bibitem{BLP} D.J. Bernstein, T. Lange, C. Peters, Attacking and Defending the McEliece Cryptosystem, in: J. Buchmann, J. Ding  (Eds.), Post-Quantum Cryptography, PQCrypto 2008, LNCS, vol. 5299.,  Springer, Berlin, Heidelberg, 2008, pp. 31–46.

\bibitem{boer0} M.A. de Boer, MDS Codes from Hyperelliptic Curves. in: R. Pellikaan, M. Perret, S.G. Vl{\"a}dut  (Eds.), Arithmetic, Geometry and Coding Theory, Walter de Gruyter, Berlin, New York, 1996, pp. 23-34. 

\bibitem{boer1} M.A. de Boer, Codes: their parameters and geometry, Eindhoven, Technische Universiteit Eindhoven, 1997, DOI: 10.6100/IR492527. 

\bibitem{boer} M.A. de Boer, The generalized Hamming weights of some hyperelliptic codes, J. Pure Appl. Algebra 123 (1998) 153-163.  

\bibitem{BrillNoether} A. von Brill,  M. Noether, \"Uber die algebraischen Functionen und ihre Anwendung in der Geometrie, Math. Annalen 7 (1874) 269-316.  


\bibitem{Cantor} D.G. Cantor, Computing in the Jacobian of a Hyperelliptic Curve, 
Math. Comp. 48 (1987) 95-101.  

\bibitem{CastellanosFanali} A.S. Castellanos, G.C. Tizziotti, Two-Point AG Codes on the GK Maximal Curves, 
IEEE Trans. Inf. Theory 62 (2016) 681-686.  

\bibitem{DinhMooreRussell} H. Dinh, C. Moore, A. Russell, McEliece and Niederreiter Cryptosystems That Resist Quantum Fourier Sampling Attacks, in: P. Rogaway (Ed.), Advances in Cryptology – CRYPTO 2011, CRYPTO 2011, LNCS, vol. 6841., Springer, Berlin, Heidelberg, 2011, pp. 761-779. 

\bibitem{FanaliGiulietti} S. Fanali, M. Giulietti, One-Point AG Codes on the GK Maximal Curves, IEEE
Trans. Inf. Theory 56 (2010) 202-210. 


\bibitem{G2} V.D. Goppa, Algebraic-geometric codes, Izv. Akad. Nauk SSSR Ser. Mat. 46 (1982) 
 762-781. (in Russian) 

\bibitem{Hess}  F. Hess, Computing Riemann-Roch Spaces in Algebraic Function Fields and Related Topics, 
J. Symbolic Comp. 33 (2002) 425--445.  

\bibitem{HuangIerardi} M. Huang, D. Ierardi, Efficient algorithms for the Riemann-Roch problem and for addition in the Jacobian of a curve, J. Symbolic Comp. 18 (1994) 519-539.  

\bibitem{JanwaMoreno} H. Janwa, O. Moreno, McEliece public key cryptosystems using algebraic-geometric codes, Des. Codes Crypt. 8  (1996) 293–307.  

\bibitem{koblitz} N. Koblitz, Hyperelliptic Cryptosystems,  J. Cryptology 1 (1989) 139–150. 

\bibitem{KorchmarosNagyTimpanella} G. Korchm\'aros, G.P. Nagy, M. Timpanella, Codes and Gap Sequences of Hermitian Curves,  IEEE Trans. Inf. Theory 66 (2020) 3547--3554.   

\bibitem{KorchmarosSpeziali} G. Korchm\'aros, P. Speziali, Hermitian codes with automorphism group isomorphic to
$\mathrm{PGL}(2, q)$ with $q$ odd, Finite Fields Their Appl. 44 (2017) 1-17.  

\bibitem{kuroki} J. Kuroki, M. Gonda, K. Matsuo, J. Chao, S. Tsujii, Fast Genus Three Hyperelliptic Curve 
Cryptosystems, In The 2002 Symposium on Cryptography and Information Security, Japan - SCIS 2002, 2002.

\bibitem{Lange2}  T. Lange, Formulae for Arithmetic on Genus 2 Hyperelliptic Curves, AAECC 15 (2005) 295-328.  

\bibitem{brigand} D. Le Brigand, Decoding of codes on hyperelliptic curves, in: G. Cohen,  P. Charpin (Eds.), EUROCODE '90, EUROCODE 1990, LNCS, vol. 514., Springer, Berlin, Heidelberg, 1990, pp. 126--134. 

\bibitem{GluherSpaenlehauer} A. Le Gluher,  P.-J. Spaenlehauer, A fast randomized geometric algorithm for computing Riemann-Roch spaces, Math. Comp. 89 (2020) 2399-2433.  

\bibitem{Lockhart} P. Lockhart, On the discriminant of a hyperelliptic curve, Trans. Amer. Math. Soc. 342 (1994) 729-752.  

\bibitem{MMPR} I. Márquez-Corbella, E. Martínez-Moro, R. Pellikaan, D. Ruano, Computational aspects of retrieving a representation of an algebraic geometry code, J. Symbolic Comp. 64 (2014) 67-87.  

\bibitem{matsuo} K. Matsuo, J. Chao, S. Tsujii, Fast Genus Two Hyperelliptic Curve Cryptosystems, ISEC 2001-31, IEICE 2001.

\bibitem{Mc} R.J. McEliece, A Public-Key Cryptosystem Based On Algebraic Coding Theory, DSN Progress Report 44 (1978) 114-116.  

\bibitem{niehage2} A. Niehage, Nonbinary Quantum Goppa Codes Exceeding the Quantum Gilbert-Varshamov Bound, Quantum Inf. Process 6 (2007) 143-158.  

\bibitem{pelzl} J. Pelzl, T. Wollinger, J. Guajardo, C. Paar, Hyperelliptic curves cryptosystems: closing the performance gap to elliptic curves,  in: C.D. Walter, Ç.K. Koç, C. Paar  (Eds.) Cryptographic Hardware and Embedded Systems - CHES 2003, CHES 2003, LNCS, vol. 2779., Springer, Berlin, Heidelberg, 2003, pp. 351-365. 

 
\bibitem{stepanov} S.A. Stepanov, Codes on fibre products of hyperelliptic curves, 
Disc. Math. Appl. 7 (1997) 77--88.  

\bibitem{Stichtenoth} H. Stichtenoth, Algebraic function fields and codes, Springer, Berlin, Heidelberg, 2009. 

\bibitem{sutherland} A.V. Sutherland, Fast Jacobian arithmetic for hyperelliptic curves of genus $3$, Thirteenth Algorithmic Number Theory Symposium (ANTS XIII), Open Book Series 2, 2019, pp. 425-442.

\bibitem{yaghoobian} T. Yaghoobian, I.F. Blake, Codes from hyperelliptic curves, in: Proc. 30-th Allerton Conf. Comm., Control, and Computing, Monticello, IL., October 1992.




\end{thebibliography}
\end{document}